\documentstyle[12pt, amsfonts]{amsart}
\begin{document}
\def\R{{\mathbb R}}
\def\Z{{\mathbb Z}}
\def\C{{\mathbb C}}
\newcommand{\trace}{\rm trace}
\newcommand{\Ex}{{\mathbb{E}}}
\newcommand{\Prob}{{\mathbb{P}}}
\newcommand{\E}{{\cal E}}
\newcommand{\F}{{\cal F}}
\newtheorem{df}{Definition}
\newtheorem{theorem}{Theorem}
\newtheorem{lemma}{Lemma}
\newtheorem{pr}{Proposition}
\newtheorem{co}{Corollary}
\def\n{\nu}
\def\sign{\mbox{ sign }}
\def\a{\alpha}
\def\N{{\mathbb N}}
\def\A{{\cal A}}
\def\L{{\cal L}}
\def\X{{\cal X}}
\def\F{{\cal F}}
\def\c{\bar{c}}
\def\v{\nu}
\def\d{\delta}
\def\diam{\mbox{\rm dim}}
\def\vol{\mbox{\rm Vol}}
\def\b{\beta}
\def\t{\theta}
\def\l{\lambda}
\def\e{\varepsilon}
\def\colon{{:}\;}
\def\pf{\noindent {\bf Proof :  \  }}
\def\endpf{ \begin{flushright}
$ \Box $ \\
\end{flushright}}

\title[Hyperplane inequality for measures]
{A hyperplane inequality for measures of unconditional convex bodies}

\author{Alexander Koldobsky}

\address{Department of Mathematics\\ 
University of Missouri\\
Columbia, MO 65211}

\email{koldobskiya@@missouri.edu}

\begin{abstract}  
We prove that there exists an absolute constant $C$ such that for every $n\in \N,$
every unconditional convex body $L$ in $\R^n$
and every measure $\mu$ with non-negative even continuous
density in $\R^n,$
$$\mu(L)\ \le\ C \max_{\xi \in S^{n-1}} 
\mu(L\cap \xi^\bot)\ |L|^{1/n} \ ,$$
where  $\xi^\bot$ is the central hyperplane in $\R^n$ perpendicular to $\xi,$ and
$|L|$ is the volume of $L.$ This is an extension to arbitrary measures of the hyperplane 
inequality for volume of unconditional convex bodies originally proved by Bourgain.
The proof is based on stability inequalities for intersection bodies. 
We also prove a similar inequality for duals of bodies with bounded
volume ratio.
\end{abstract}  
\maketitle

\section{Introduction}
The hyperplane problem \cite{Bo1, Bo2, Ba1, MP}, a major open problem in convex geometry,
asks whether there exists an absolute constant $C$ so that for any origin-symmetric convex body $K$ in $\R^n$
of volume 1 there is a hyperplane section of $K$ whose $(n-1)$-dimensional volume is greater than $1/C.$
In other words, does there exist a constant $C$ so that for any $n\in \N$ and any
origin-symmetric convex body $K$ in $\R^n$
\begin{equation} \label{hyper}
|K|^{\frac {n-1}n} \le C \max_{\xi \in S^{n-1}} |K\cap \xi^\bot|,
\end{equation}
where  $\xi^\bot$ is the central hyperplane in $\R^n$ perpendicular to $\xi,$ and
$|K|$ stands for volume of proper dimension?
The best current result $C\le O(n^{1/4})$ is due to Klartag \cite{Kl}, who
removed the  logarithmic term from an earlier estimate of Bourgain \cite{Bo3}.
We refer the reader to [BGVV] for the history and partial results.

For certain classes of bodies the question has been answered in affirmative. These classes
include unconditional convex bodies (as initially observed by Bourgain; see \cite{MP, J2,
BN, BGVV} for different proofs), unit balls of subspaces of $L_p$  \cite{Ba2, J1, M}, intersection bodies
\cite[Th. 9.4.11]{G}, zonoids, duals of bodies with bounded volume ratio
\cite{MP}, the Schatten classes \cite{KMP}, $k$-intersection bodies \cite{KPY, K8}.

In this note we prove that volume in (\ref{hyper}) can be replaced by an arbitrary measure 
for two classes of bodies - unconditional convex bodies and duals of convex bodies with bounded
volume ratio. Namely, there exists an absolute constant $C$ such that for any body $K$ from these 
classes and any measure $\mu$ with even continuous density $f$ in $\R^n$ 
\begin{equation} \label{arbmeasure}
\mu(K) \le C \max_{\xi \in S^{n-1}} \mu(K\cap \xi^\bot)\ |K|^{1/n},
\end{equation}
where $\mu(B)=\int_B f$ for every compact set $B\subset \R^n.$

For intersection bodies inequality (\ref{arbmeasure}) was proved in \cite{K4} with the best possible constant  
$C=\frac n{n-1} c_n,$ where $c_n= |B_2^n|^{\frac {n-1}n} / |B_2^{n-1}|<1$ and $B_2^n$ is the unit Euclidean ball in $\R^n.$
For $k$-intersection bodies, (\ref{arbmeasure}) was proved in \cite{K8} with $C$ depending only on $k.$ 
For general origin-symmetric convex bodies,
it is possible to prove (\ref{arbmeasure}) with $C\le O(\sqrt{n})$ (see \cite{K5}):
\begin{equation} \label{sqrtn2}
\mu(K)\ \le\ \sqrt{n} \frac n{n-1} c_n\max_{\xi \in S^{n-1}} 
\mu(K\cap \xi^\bot)\ |K|^{1/n}.
\end{equation}
Analogs of (\ref{sqrtn2}) for sections of lower dimensions and for complex convex bodies 
were proved in \cite{K6}. It was shown in \cite{K8}  that the constant $\sqrt{n}$ in (\ref{sqrtn2}) 
can be replaced by $n^{1/2-1/p}$ when
$K$ is the unit ball of an $n$-dimensional subspace of $L_p,\ p>2.$ The author does not know whether the 
constant in (\ref{sqrtn2}) is optimal for arbitrary measures; 
for volume it is certainly not optimal.
 
The proofs are based on a stability result for intersection bodies \cite{K4}. This result
continues the study of stability in volume comparison problems initiated in \cite{K3, K7}.

\section{Unconditional convex bodies}

We need several definitions and facts.
A closed bounded set $K$ in $\R^n$ is called a {\it star body}  if 
every straight line passing through the origin crosses the boundary of $K$ 
at exactly two points different from the origin, the origin is an interior point of $K,$
and the {\it Minkowski functional} 
of $K$ defined by 
$$\|x\|_K = \min\{a\ge 0:\ x\in aK\}$$
is a continuous function on $\R^n.$ 

The {\it radial function} of a star body $K$ is defined by
$$\rho_K(x) = \|x\|_K^{-1}, \qquad x\in \R^n,\ x\neq 0.$$
If $x\in S^{n-1}$ then $\rho_K(x)$ is the radius of $K$ in the
direction of $x.$

The class of intersection bodies was introduced by Lutwak \cite{L}.
Let $K, L$ be origin-symmetric star bodies in $\R^n.$ We say that $K$ is the 
intersection body of $L$ if the radius of $K$ in every direction is 
equal to the $(n-1)$-dimensional volume of the section of $L$ by the central
hyperplane orthogonal to this direction, i.e. for every $\xi\in S^{n-1},$
\begin{equation} \label{intbodyofstar}
\rho_K(\xi)= \|\xi\|_K^{-1} = |L\cap \xi^\bot|.
\end{equation} 
All bodies $K$ that appear as intersection bodies of different star bodies
form {\it the class of intersection bodies of star bodies}. The class of {\it intersection bodies} 
can be defined as the closure of the class of intersection bodies of star bodies
in the radial metric 
$$\rho(K,L)=\sup_{\xi\in S^{n-1}} \left|\rho_K(\xi)-\rho_L(\xi)\right|.$$  
Intersection bodies played a crucial role in the solution of the
Busemann-Petty problem and its generalizations; see \cite[Ch. 5]{K1}.

We use  the following stability result for intersection bodies proved in \cite{K5}.  
\begin{pr}\label{stab1}{\bf (\cite{K5})}
Suppose that $K$ is an intersection body in $\R^n,$  $f$
is an even continuous function on $K,$ $f\ge 1$ everywhere on $K,$ and $\e>0.$ If
$$
\int_{K\cap \xi^\bot} f \ \le\ |K\cap \xi^\bot| +\e,\qquad \forall \xi\in S^{n-1},
$$
then
$$
\int_K f\ \le\ |K| + \frac {n}{n-1}\ c_{n}\ |K|^{1/n}\e;
$$
recall that $c_n<1.$
\end{pr}

Let $e_i,\ 1\le i\le n,$ be the standard basis of $\R^n.$ A star body $K$ in $\R^n$ is called unconditional 
if for every choice of real numbers $x_i$ and $\delta_i=\pm 1,\ 1\le i \le n$
we have 
$$\|\sum_{i=1}^n \delta_ix_i e_i \|_K = \|\sum_{i=1}^n x_i e_i \|_K.$$

\begin{theorem}\label{uncond} There exists an absolute constant $C$ such that for every $n\in \N,$
every unconditional convex body $L$ in $\R^n$
and every measure $\mu$ with even continuous non-negative density on $L$
\begin{equation} \label{uncond1}
\mu(L)\ \le\  C \max_{\xi \in S^{n-1}}  \mu(L\cap \xi^\bot)\ |L|^{1/n}.
\end{equation}
\end{theorem}

\pf  By a result of Lozanovskii \cite{Lo} (see the proof in \cite[Corollary 3.4]{P}), there exists a (diagonal) linear
operator $T: \R^n\to \R^n$ so that 
$$T(B_\infty^n) \subset L \subset n T(B_1^n),$$ 
where $B_p^n$ is the unit ball of the space $\ell_p^n.$
Let $K=nT(B_1^n).$ By \cite[Th. 3]{K2} and the fact that a linear transformation of an intersection body
is an intersection body (see \cite{L}; it also follows from \cite[Th. 1]{K2} and the connection between linear
transformations and the Fourier transform), the body $K$ is an intersection body in $\R^n.$
Let $g$ be the density of $\mu,$ and let $f= \chi_K + g \chi_L,$ where $\chi_K,\ \chi_L$ are the indicator functions of $K$ and $L.$ 
Then $f\ge 1$ everywhere on $K.$ Put 
$$\e=\max_{\xi\in S^{n-1}} \left(\int_{K\cap \xi^\bot} f - |K\cap \xi^\bot| \right)= \max_{\xi\in S^{n-1}} \int_{L\cap \xi^\bot} g.$$
Now we can apply Proposition \ref{stab1} to $f,K,\e$ (the function $f$ is not necessarily continuous on $K,$ 
but the result holds by a simple approximation argument). We get
$$\mu(L)= \int_L g = \int_K f  -\ |K|$$
$$ \le \frac n{n-1} c_{n}  |K|^{1/n}\max_{\xi\in S^{n-1}} \int_{L\cap \xi^\bot} g= 
\frac n{n-1} c_{n}  |K|^{1/n}\max_{\xi\in S^{n-1}} \mu(L\cap \xi^\bot) .$$
Since $|B_1^n| = 2^n/n!$ (see for example \cite[Lemma 2.19]{K1}), we have
$|K|^{1/n}\le |\det T|^{1/n} e.$ On the other hand, $|T(B_\infty^n)| = 2^n |\det T|,$ and $T(B_\infty^n)\subset L,$
so $|K|^{1/n} \le \frac e2\  |L|^{1/n}.$ Since $c_n<1,$ (\ref{uncond1}) follows with $C=e$ (roughly estimating 
the constant; in fact the constants that we get tend to $\frac 12 \sqrt{e}$ as $n\to \infty)$. \qed 

\section{Duals of bodies with bounded volume ratio}

The volume ratio of a convex body $K$ in $\R^n$ is defined by
$$\rm{v.r.}(K) =\inf_E \left\{ \left(\frac{|K|}{|E|}\right)^{1/n}:\ E\subset K,\ E-{\rm ellipsoid}\right\}.$$
The following argument is standard and first appeared in \cite{BM} and \cite{MP}.
Let $K^\circ$ and $E^\circ$ be polar bodies of $K$ and $E,$ respectively. If $E$ is an ellipsoid,
then 
$$|E| |E^\circ|=|B_2^n|^2.$$
By the reverse Santalo inequality of Bourgain and Milman \cite{BM}, there exists an absolute constant
$c>0$ such that
$$\left(|K||K^\circ|\right)^{1/n} \ge \frac cn.$$
Combining these and using the asymptotics of $B_2^n$ we get that there exists an absolute
constant $C$ such that
$$ \left(\frac{|E^\circ|}{|K^\circ|}\right)^{1/n} \le C \left(\frac{|K|}{|E|}\right)^{1/n}.$$

\begin{theorem}\label{boundedvr} There exists an absolute constant $C$ such that for every $n\in \N,$
every origin-symmetric convex body $L$ in $\R^n$
and every measure $\mu$ with even continuous non-negative density on $L$
\begin{equation} \label{uncond1}
\mu(L)\ \le\  C\ \rm{v.r.}(L^\circ) \max_{\xi \in S^{n-1}}  \mu(L \cap \xi^\bot)\ |L|^{1/n}.
\end{equation}
\end{theorem}

\pf If $E$ is an ellipsoid, $E\subset L^\circ$, then the ellipsoid $E^\circ$ contains $L.$
Also every ellipsoid  is an intersection body as a linear image of the Euclidean ball.
Applying Proposition \ref{stab1} in the same way as in Theorem \ref{uncond} (with $K=E^\circ)$
and using the argument before the statement of the theorem, we get
$$\mu(L)\le  \frac n{n-1} c_{n}  \max_{\xi\in S^{n-1}} \mu(L\cap \xi^\bot)|E^\circ|^{1/n}$$
$$\le C   \frac n{n-1} c_{n} \max_{\xi\in S^{n-1}} \mu(L\cap \xi^\bot)\left(\frac{|L^\circ|}{|E|}\right)^{1/n} |L|^{1/n}.$$
The result follows. \qed

\bigbreak
{\bf Acknowledgement.} I wish to thank the US National Science Foundation for support through 
grant DMS-1265155.

\end{document}